\documentclass[11pt,reqno]{amsart}
\usepackage[a4paper,margin=28mm]{geometry}
\usepackage[T1]{fontenc}
\usepackage{lmodern,microtype,amsmath,amssymb,amsthm,mathtools}
\usepackage[hidelinks]{hyperref}
\hypersetup{pdftitle={Infinite Generalized Sidon Sets with a Negative Answer to a Problem of Szemeredi},pdfauthor={Guoqing Wang and Wenkai Yang}}
\newtheorem{theorem}{Theorem}[section]
\newtheorem{lemma}[theorem]{Lemma}
\newtheorem{proposition}[theorem]{Proposition}
\newtheorem{corollary}[theorem]{Corollary}
\theoremstyle{remark}

\newcommand{\N}{\mathbb N}
\newcommand{\F}{\mathbb F}
\newcommand{\Z}{\mathbb Z}

\numberwithin{equation}{section}
\title[Infinite generalized Sidon sets]{Infinite Generalized Sidon Sets with\\
a Negative Answer to a Problem of Szemer\'edi}
\author{Guoqing Wang \and Wenkai Yang}
\date{10 November 2026}
\address{School of Mathematical Sciences\\Tiangong University\\Tianjin 300387, P.R. China}
\email{gqwang1979@aliyun.com}
\address{Center for Combinatorics\\Nankai University\\Tianjin 300071, P.R. China}
\email{9820250116@nankai.edu.cn}


\begin{document}
\begin{abstract}
	Let $A(x)=|A\cap[1,x]|$.
	For every integer $h\ge3$, every finite $B_h$ set
	$A_0\subseteq\N$, and every $0<c\le h^{-2/h}$,
	we construct an infinite $B_h$ set $A\subseteq\N$ containing
	$A_0$ such that $A(x)/x^{1/h}$ has upper limit $c$
	and lower limit zero.
	This extends the range $0<c\le1/h$ obtainable from our first
	construction, which answered a question of Cilleruelo and
	Nathanson, later raised by O'Bryant, on the existence of a
	positive upper limit.
	For $h=3$, a positive upper limit disproves the assertion
	$A(x)=o(x^{1/3})$ in Szemer\'edi's Problem~186(d).
	The improved range follows from an integer lift of a modular
	construction of G\'omez Ruiz and Trujillo Solarte, with large
	gaps between sums of distinct $r$-element multisets for each
	fixed $1\le r<h$.
	These gaps ensure that successive blocks can be adjoined
	by translation while preserving the $B_h$ property.
	The same prescribed limits can be attained inside every set
	of asymptotic density one containing $A_0$.
	Separately, we construct $2^{\aleph_0}$ extensions of $A_0$
	that are maximal among all $B_h$ subsets of $\N$,
	have the prescribed limits, and avoid any prescribed finite
	subset of $\N\setminus A_0$.
	These maximal extensions can be chosen so that
	$\sum_{a\in A}a^{-1/h}$ diverges arbitrarily slowly.
\end{abstract}
\maketitle

\section{Introduction}

Throughout the paper, $h$ is an integer with $h\ge2$.
Let $\N=\{1,2,\ldots\}$. A set $A\subseteq\Z$ is a $B_h$ set if
$a_1+\cdots+a_h=a'_1+\cdots+a'_h$, with all summands in $A$, implies
equality of the two multisets of summands. Repetitions are permitted.
For $A\subseteq\N$, let $A(x)=|A\cap[1,x]|$. A nonempty $B_h$ set
is also a $B_r$ set for $1\le r<h$, by adding the same $h-r$ terms to both
sides of an equality. The distinct $h$-term sums of a finite
$B_h$ set $S\subseteq[1,n]$ give
\begin{equation}\label{eq:counting}
 \frac{|S|^h}{h!}\le\binom{|S|+h-1}{h}\le hn,
 \qquad |S|\le(h!hn)^{1/h}.
\end{equation}
Thus the natural density scale is $x^{1/h}$.

Finite constructions of Bose and Chowla~\cite{bose-chowla} attain this
order of magnitude. The infinite problem requires dense finite sets at
different scales to be compatible. Kr\"uckeberg~\cite{krueckeberg}
constructed an infinite Sidon set with normalized upper limit at least
$1/\sqrt2$.
Cilleruelo and Nathanson~\cite{cn-representation} observed that, for
each $h\ge3$ and each function $\omega(x)\to\infty$, one can construct
an infinite $B_h$ set $G$ with
$\limsup_{x\to\infty}\omega(x)G(x)/x^{1/h}>0$ by dilating successive blocks.
Their 2007 preprint identified the existence of a positive upper limit
without the factor $\omega(x)$ as open. O'Bryant~\cite[v1, Section~5]{obryant2026}
later conjectured that such a set exists for every $h\ge3$.
In the first version of this paper,
we constructed an extension of every finite $B_h$ set with normalized upper
limit $1/h$ and lower limit zero. This settled the existence question
for $h\ge3$, as O'Bryant records in~\cite[v2, Section~5.3]{obryant2026}.
Proportional thinning while retaining the initial set gives each
$0<c\le 1/h$. Here we enlarge that range to $0<c\le h^{-2/h}$ for $h\ge3$.

We include $h=2$ in Theorem~\ref{thm:main} for use in the extensions below.
In that case the statement follows from earlier work.
O'Bryant's proof of Kr\"uckeberg's bound can start with any finite Sidon set
by~\cite[Lemma~9]{obryant-I}. The normalized lower limit is zero for every
Sidon set by Erd\H{o}s's bound (see~\cite[Theorem~1]{obryant-I}).
If the upper limit is $b\ge1/\sqrt2$, retaining the initial set and a
proportion $c/b$ of the remaining terms gives each $0<c\le1/\sqrt2$.

\begin{theorem}\label{thm:main}
Let $h\ge2$, let $A_0\subseteq\N$ be a finite $B_h$ set, and define
\[
 c_h=\begin{cases}1/\sqrt2,&h=2,\\h^{-2/h},&h\ge3.\end{cases}
\]
For every $0<c\le c_h$, there is an infinite $B_h$ set
$A\subseteq\N$ containing $A_0$ such that
\begin{equation}\label{eq:limits}
 \limsup_{x\to\infty}\frac{A(x)}{x^{1/h}}=c,
 \qquad \liminf_{x\to\infty}\frac{A(x)}{x^{1/h}}=0.
\end{equation}
\end{theorem}

For $h=3$, our first version disproved the assertion
$A(x)=o(x^{1/3})$ in
Szemer\'edi's Problem~186(d)~\cite[p.~56]{szemeredi}
by attaining upper limit $1/3$.
Theorem~\ref{thm:main} raises this value to $3^{-2/3}$.
This full-limit assertion is distinct from the question whether
every infinite $B_3$ set has normalized lower limit zero.
Helm~\cite{helm} discusses the latter question and rules out
$A(x)\sim\alpha x^{1/3}$ for $\alpha>0$, while
Chen~\cite{chen} gives logarithmic lower-limit restrictions
for even orders.
Neither conclusion excludes a positive normalized upper limit.

The finite initial set can also be prescribed inside a restricted
ambient set. Density below always means asymptotic density in $\N$.

\begin{theorem}\label{thm:density-one}
Let $E\subseteq\N$ have density one, and let $A_0\subseteq E$ be a
finite $B_h$ set. For every $0<c\le c_h$, there is an infinite
$B_h$ set $A$ with $A_0\subseteq A\subseteq E$ satisfying
\eqref{eq:limits}.
\end{theorem}

In particular, a prescribed set of density zero can be avoided whenever
it is disjoint from the initial set. Section~\ref{sec:ambient} also
treats positive upper density, infinite initial sets, and compatible
families. The constants are attainable values; optimality is not asserted.

Maximality under inclusion in $\N$ can also coexist with the
limits in \eqref{eq:limits}. A $B_h$ set is maximal if no further
positive integer can be adjoined while preserving the $B_h$ property.

\begin{theorem}\label{thm:maximal-intro}
	Let $A_0\subseteq\N$ be a finite $B_h$ set, and let
	$F\subseteq\N\setminus A_0$ be finite.
	For every $0<c\le c_h$, there are $2^{\aleph_0}$ infinite
	maximal $B_h$ sets $A\subseteq\N$ containing $A_0$,
	disjoint from $F$, and satisfying \eqref{eq:limits}.
\end{theorem}

Theorem~\ref{thm:maximal} proves this statement together with
all the conclusions of Proposition~\ref{prop:slow-scales}.
In contrast, Proposition~\ref{prop:sparse-maximal} constructs
arbitrarily sparse maximal sets, showing why sparsity alone
does not suffice for extending an infinite initial set.

The main issue in joining blocks is the separation of shorter sums.
A large translation distinguishes sums with different numbers of new
terms. With the same number $1\le r<h$, however, the translation
cancels, and a mixed relation becomes
\[
 \sum_{i=1}^r e_i-\sum_{i=1}^r e'_i
 =\sum_{j=1}^{h-r}s'_j-\sum_{j=1}^{h-r}s_j.
\]
If the old set lies in $[1,n]$, the right side has absolute value
at most $hn$.
For each fixed $1\le r<h$, our finite blocks ensure that sums
of distinct $r$-element multisets differ by more than $hn$. In our first version we encoded the finite-field
moment curve of~\cite{hajela,ruiz-trujillo} in base $hq$ and proved the
separation estimate that yields upper limit $1/h$. For the improved
constant $h^{-2/h}$ we use the modular logarithmic graph of G\'omez Ruiz
and Trujillo Solarte~\cite[Theorem~4]{gomez-trujillo}. Its integer lift
has shorter-sum gaps large enough for the same iteration and fits in a
shorter interval. For $h=2$ we apply the deletion argument
of Cilleruelo and Nathanson~\cite[Lemma~3.2]{cn-perfect} to a
Bose--Chowla block~\cite{bose-chowla}.

The study of Sidon sets goes back to Erd\H{o}s and Tur\'an~\cite{erdos-turan};
see~\cite{halberstam-roth,obryant-survey} for broader accounts.
Ajtai, Koml\'os and Szemer\'edi~\cite{aks}, Ruzsa~\cite{ruzsa}, and
Cilleruelo~\cite{cilleruelo} obtain lower bounds holding at every large
scale. For $h\ge3$, Cilleruelo constructs a $B_h$ set with
$A(x)=x^{\sqrt{(h-1)^2+1}-(h-1)+o(1)}$.
Our growth statements in Section~\ref{sec:growth} concern simultaneous
control of the prescribed positive upper limit, a finite initial set,
and the behavior between dense scales; they do not claim a larger
all-scale exponent than those earlier results.

Section~\ref{sec:blocks} gives one recursive construction and proves
the three available finite-block choices and their compatibility.
Section~\ref{sec:density} proves Theorem~\ref{thm:main} and the
exact value $1/h$ for the moment-curve choice.
Sections~\ref{sec:growth} and~\ref{sec:ambient} apply the common
lemmas to growth, ambient restrictions and prescribed initial sets.
Section~\ref{sec:maximal} treats maximality, the number of extensions,
and the obstruction from sparse maximal initial sets.

\section{Construction and finite blocks}\label{sec:blocks}

We give one construction for the three block choices, verify it below,
and prove its density in Section~\ref{sec:density}.
For an integer $a$, let $[a]_q\in\{0,\ldots,q-1\}$ be its
representative modulo $q$.

\medskip
\noindent\textbf{Step 1. The old set and the prime.}
Start with a finite nonempty $B_h$ set $S\subseteq[1,n]$, where
$n$ is a positive integer, and choose a prime $q>hn$.
The value $n=\max S$ is permitted. The parameters satisfy
\begin{equation}\label{eq:parameter-order}
 h\le hn<q,\qquad
 \frac{|S|}{q}\le\frac{(h!hn)^{1/h}}q
 <(h!)^{1/h}q^{1/h-1}.
\end{equation}

\medskip
\noindent\textbf{Step 2. The finite block.}
Choose a nonempty block $E_q\subseteq[0,dq^h)$ of integers,
where $d\ge1$ is an integer fixed throughout the recursion.
The moment-curve choice is
\begin{equation}\label{eq:moment-block}
 f_q(a)=\sum_{j=1}^h[a^j]_q(hq)^{h-j},\qquad
 E_q=\{f_q(a):0\le a<q\},\qquad d=h^{h-1}.
\end{equation}
The other choices in Lemmas~\ref{lem:logarithmic}
and~\ref{lem:sidon} have $d=h$ and $d=1$, respectively;
the second is available only for $h=2$.

\medskip
\noindent\textbf{Step 3. The required block properties.}
Each chosen $E_q$ is a $B_h$ set, and sums of distinct
$r$-element multisets from it differ by more than $hn$ whenever
$1\le r<h$. For the block in \eqref{eq:moment-block},
$|E_q|=q$, and for distinct parameter multisets,
\begin{equation}\label{eq:moment-gap}
 \left|\sum_{i=1}^r f_q(a_i)-\sum_{i=1}^r f_q(a'_i)\right|
 \ge hq-r(q-1)>q>hn\qquad(1\le r<h).
\end{equation}

\medskip
\noindent\textbf{Step 4. The translation and the enlarged set.}
Define
\begin{equation}\label{eq:one-step}
 M=(h-1)dq^h+q,\qquad S'=S\cup(M+E_q).
\end{equation}
Then $S'$ is a $B_h$ set,
$n<M\le M+e<hdq^h+q$ for every $e\in E_q$, and
$|S'|=|S|+|E_q|$.
Its normalized count at the endpoint is
\[
 \frac{|S'|}{(hdq^h+q)^{1/h}}
 =\frac{|S|/q+|E_q|/q}{(hd+q^{1-h})^{1/h}}.
\]

\medskip
\noindent\textbf{Step 5. Iteration.}
Start with $S_0=A_0$ if $A_0\ne\varnothing$, and with $\{1\}$
otherwise. At stage $k$, choose $q_k>h\max S_{k-1}$ and a block
satisfying Step~3 with $n=\max S_{k-1}$, and define
\begin{equation}\label{eq:recursion}
 \begin{aligned}
 M_k&=(h-1)dq_k^h+q_k,& n_k&=hdq_k^h+q_k,\\
 S_k&=S_{k-1}\cup(M_k+E_{q_k}),&
 U&=\bigcup_{k\ge0}S_k.
 \end{aligned}
\end{equation}
Here $n_k$ is an upper endpoint, not necessarily an element of $S_k$.
Larger bounds for the old set and correspondingly larger primes may
also be used. The set $U$ is an infinite $B_h$ set containing $A_0$.

\medskip
\noindent\textbf{Step 6. Density.}
Whenever $|E_{q_k}|/q_k\to1$, the source set satisfies
\begin{equation}\label{eq:source-scales}
 \frac{U(n_k)}{n_k^{1/h}}\longrightarrow(hd)^{-1/h},\qquad
 \frac{U(M_k-1)}{(M_k-1)^{1/h}}\longrightarrow0,\qquad
 \frac{n_k}{M_k-1}\longrightarrow\frac h{h-1}.
\end{equation}
For the moment-curve choice the normalized upper limit is exactly
$1/h$. Using the logarithmic block for $h\ge3$ and the Sidon block
for $h=2$, selection gives every prescribed upper limit up to $c_h$.
In the Sidon deletion construction,
$|E_{q_k}|/q_k\to1$ is ensured by making the old bound negligible
relative to $q_k$.

\medskip
The next lemma proves that the extension in \eqref{eq:one-step}
is a $B_h$ set and also permits the larger translations used later
inside a prescribed ambient set.

\begin{lemma}\label{lem:extension}
Let $S\subseteq[1,n]$ be a finite nonempty $B_h$ set, and let
$E\subseteq[0,D)$ be a finite nonempty $B_h$ set of integers,
where $n,D$ are positive integers.
Suppose that sums of distinct $r$-element multisets from $E$ differ
by more than $hn$ for every $1\le r<h$.
For every integer $M\ge(h-1)D+hn+1$, the set $S\cup(M+E)$
is a $B_h$ set, and the new block lies above $S$.
\end{lemma}
\begin{proof}
We first prove that equal sums contain the same number of new terms.
Since $M>n$, the old and new sets are disjoint.
If one $h$-term sum contains $r$ new terms and the other contains
$s<r$, their difference is at least
\[
 (r-s)M-sD-(h-s)n
 \ge M-(h-1)D-hn>0.
\]
Here the first sum is at least $rM$, whereas the second is
at most $s(M+D)+(h-s)n$; also $s\le h-1$.
Exchanging the sums excludes the case $r<s$.
If $r=0$ or $r=h$, the $B_h$ property of $S$ or $E$ identifies
the summands, after cancellation of the translation when $r=h$.

We resolve the remaining case $1\le r<h$ by the shorter-sum gap.
Let $e_1,\ldots,e_r$ and $e'_1,\ldots,e'_r$ be the new parameters,
and let $s_1,\ldots,s_{h-r}$ and $s'_1,\ldots,s'_{h-r}$ be the
old terms on the two sides. Cancellation gives
\[
 \sum_{i=1}^r e_i-\sum_{i=1}^r e'_i
 =\sum_{j=1}^{h-r}s'_j-\sum_{j=1}^{h-r}s_j,\qquad
 \left|\sum_{j=1}^{h-r}s'_j-\sum_{j=1}^{h-r}s_j\right|
 \le(h-r)(n-1)\le hn.
\]
The gap hypothesis forces the new multisets to agree.
After cancelling them, add the same $r$ copies of any fixed element
of the nonempty set $S$ to both old sums.
The $B_h$ property of $S$ identifies the padded multisets;
removing the added copies identifies the original old multisets.
\end{proof}

The following lemma recalls the integer encoding and shorter-sum
estimate from Proposition~2.3 of our first version. Its power-sum argument
comes from the finite-field moment curve~\cite{hajela,ruiz-trujillo}.
For $h=2$, the integer encoding $2qa+[a^2]_q$ is the
Erd\H{o}s--Tur\'an construction~\cite{erdos-turan}.

\begin{lemma}\label{lem:moment}
For a prime $q>h$, the set in \eqref{eq:moment-block} is a $B_h$
set of $q$ elements, the map $a\mapsto f_q(a)$ is strictly increasing,
and
\begin{equation}\label{eq:moment-range}
 a(hq)^{h-1}\le f_q(a)<(a+1)(hq)^{h-1}
 \le h^{h-1}q^h\qquad(0\le a<q).
\end{equation}
For $1\le r<h$, sums of distinct $r$-element parameter multisets
differ by at least $hq-r(q-1)>q$ after applying $f_q$.
\end{lemma}
\begin{proof}
We first prove the range and cardinality assertions.
Since $[a]_q=a$ for $0\le a<q$, the terms after
$a(hq)^{h-1}$ in $f_q(a)$ satisfy
\[
 0\le f_q(a)-a(hq)^{h-1}
 \le\frac{q-1}{hq-1}\bigl((hq)^{h-1}-1\bigr)
 <(hq)^{h-1}.
\]
This proves \eqref{eq:moment-range}.
The intervals for consecutive parameters are disjoint and ordered,
so $f_q$ is strictly increasing and $|E_q|=q$.

We next identify the power sums that distinguish two parameter multisets.
For two $r$-tuples in $\F_q$ with common power sums $p_i$, $1\le i\le r<q$,
let $e_j,e'_j$ be their elementary symmetric functions and
$e_0=e'_0=1$. Newton's identities give
\[
 j(e_j-e'_j)
 =\sum_{i=1}^j(-1)^{i-1}(e_{j-i}-e'_{j-i})p_i
 \qquad(1\le j\le r).
\]
Since each $j$ is invertible in $\F_q$, induction gives $e_j=e'_j$.
The two monic polynomials with the tuple entries as roots are
therefore equal. Unique factorization identifies each root with its
multiplicity, so the two multisets agree.

We now prove uniqueness and all shorter-sum gaps in one comparison.
Fix $1\le r\le h$ and distinct multisets of parameters
$a_1,\ldots,a_r$ and $a'_1,\ldots,a'_r$ in $\{0,\ldots,q-1\}$.
For each power $1\le j\le h$, let
$u_j=\sum_{i=1}^r[a_i^j]_q-\sum_{i=1}^r[(a'_i)^j]_q$;
then $|u_j|\le r(q-1)$.
The preceding power-sum argument shows that not all
$u_1,\ldots,u_r$ are zero.
Let $\ell\le r$ be the first index with $u_\ell\ne0$.
Its coefficient has absolute value at least one, so
\[
 \begin{aligned}
 \left|\sum_{i=1}^r f_q(a_i)-\sum_{i=1}^r f_q(a'_i)\right|
 &=\left|\sum_{j=\ell}^h u_j(hq)^{h-j}\right|\\
 &\ge(hq)^{h-\ell}
       -r(q-1)\sum_{j=0}^{h-\ell-1}(hq)^j\\
 &=1+\bigl(hq-1-r(q-1)\bigr)
          \sum_{j=0}^{h-\ell-1}(hq)^j.
 \end{aligned}
\]
An empty sum here is zero.
Since $hq-1-r(q-1)\ge h-1>0$, the last expression is at least
one; taking $r=h$ proves the $B_h$ property.
If $r<h$, then $\ell<h$ and the geometric sum is at least one.
The same expression is therefore at least
$hq-r(q-1)\ge q+h-1>q$, proving \eqref{eq:moment-gap}.
\end{proof}

For the second block we lift the modular graph of G\'omez Ruiz and
Trujillo Solarte~\cite[Theorem~4]{gomez-trujillo} to the integers.
Its order-two specialization is Ruzsa's construction~\cite{ruzsa-modular}, as observed
in~\cite[Observation~4]{gomez-trujillo}. The next lemma gives the
width of the lift and the gaps between its shorter sums.

\begin{lemma}\label{lem:logarithmic}
For a prime $q>h$, there is a $B_h$ set $E_q\subseteq[0,hq^h)$
with $q$ elements if $h\ge3$ and $q-1$ elements if $h=2$.
For $1\le r<h$, sums of distinct $r$-element multisets from $E_q$
differ by at least $hq-r(q-1)>q$.
\end{lemma}
\begin{proof}
We first define the parameters and their integer representatives.
Let $\theta$ generate the multiplicative group
$\F_{q^{h-1}}^\times$. Its degree over $\F_q$ is $h-1$:
if this degree were $\ell<h-1$, the field $\F_q(\theta)$ would
have $q^\ell$ elements, so the order $q^{h-1}-1$ of $\theta$
would be at most $q^\ell-1$, a contradiction.
For $a\in\F_q$ with $\theta+a\ne0$, let $[a]_q$ be its representative
in $\{0,\ldots,q-1\}$ and define
\begin{equation}\label{eq:logarithmic-block}
 \begin{gathered}
 \theta^{b_q(a)}=\theta+a,\qquad 0\le b_q(a)<q^{h-1}-1,\\
 f_q(a)=hq\,b_q(a)+[a]_q,\qquad
 E_q=\{f_q(a):a\in\F_q,\ \theta+a\ne0\}.
 \end{gathered}
\end{equation}
The chosen generator gives a unique integer $b_q(a)$ in the indicated
range. When $h\ge3$, the degree of $\theta$ is at least two, so
$\theta\notin\F_q$ and every $a\in\F_q$ is permitted. When $h=2$,
exactly the parameter $a=-\theta$ is excluded.

We next check the size and range in the statement.
If $f_q(a)=f_q(a')$, reduction modulo $hq$ gives $[a]_q=[a']_q$,
since both remainders lie in $[0,q)$. Thus $a=a'$.
This proves the asserted cardinality. Moreover,
\[
 0\le f_q(a)\le hq(q^{h-1}-2)+(q-1)<hq^h.
\]
These are integer inequalities, whereas the defining equation for
$\theta^{b_q(a)}$ is an equation in $\F_{q^{h-1}}$.

We prove the $B_h$ property by comparing two sums of $h$ block elements.
Let $a_1,\ldots,a_h$ and $a'_1,\ldots,a'_h$ be their respective
parameters, with repetitions allowed, and suppose that
$\sum_{i=1}^h f_q(a_i)=\sum_{i=1}^h f_q(a'_i)$.
The definition in \eqref{eq:logarithmic-block} gives
\[
 hq\left(\sum_{i=1}^h b_q(a_i)-\sum_{i=1}^h b_q(a'_i)\right)
 =-\sum_{i=1}^h[a_i]_q+\sum_{i=1}^h[a'_i]_q.
\]
The absolute value of the right side is at most $h(q-1)<hq$.
It is a multiple of $hq$, so it is zero. Both the sums of the
integer representatives $[a_i]_q$ and the sums of the integer
logarithms $b_q(a_i)$ therefore agree.

We now recover the parameters from these two equalities.
Reduction of the first equality modulo $q$ gives
$\sum a_i=\sum a'_i$ in $\F_q$, and exponentiation of $\theta$
in the second gives
$\prod(\theta+a_i)=\prod(\theta+a'_i)$.
Consequently the polynomial
\[
 \prod_{i=1}^h(X+a_i)-\prod_{i=1}^h(X+a'_i)\in\F_q[X]
\]
vanishes at $\theta$. Its coefficients of $X^h$ cancel because the
products are monic, and its coefficients of $X^{h-1}$ cancel because
the parameter sums agree. Its degree is at most $h-2$, below the
degree $h-1$ of $\theta$, so it is the zero polynomial.
Unique factorization identifies the parameter multisets, including
their multiplicities. The corresponding multisets of $f_q(a_i)$
are equal, proving the $B_h$ property. This argument also covers
$h=2$, when the difference is a constant polynomial.

We finally separate all shorter sums. Fix $1\le r<h$ and two
distinct multisets of $r$ permitted parameters. Equality of their
integer logarithm sums would imply equality of the products at
$\theta$. The difference of the two monic degree-$r$ products has
degree at most $r-1<h-1$, so the preceding degree argument would
identify the two parameter multisets, a contradiction. Thus their
logarithm sums are distinct integers and differ by at least one.
Applying the reverse triangle inequality to \eqref{eq:logarithmic-block}
gives
\[
 \left|\sum_{i=1}^r f_q(a_i)-\sum_{i=1}^r f_q(a'_i)\right|
 \ge hq-r(q-1)
 \ge hq-(h-1)(q-1)=q+h-1>q.
\]
This proves the separation required for every $r<h$.
\end{proof}

For $h=2$, the preceding block gives an attainable upper limit of $1/2$
after the selection in Lemma~\ref{lem:selection}.
The next lemma applies the short-difference deletion argument of
Cilleruelo and Nathanson~\cite[Lemma~3.2]{cn-perfect} to a
Bose--Chowla block~\cite{bose-chowla}, with the gap chosen for the old set.

\begin{lemma}\label{lem:sidon}
For a prime $q$ and a positive integer $n$, there is a Sidon set
$E_q\subseteq[0,q^2)$ with $q-2n\le|E_q|\le q$, whose distinct
elements differ by more than $2n$.
\end{lemma}
\begin{proof}
We first construct a Sidon set with $q$ elements.
Let $\theta$ generate $\F_{q^2}^\times$, and take the integer
representatives $b_q(a)$ defined by
\[
 \theta^{b_q(a)}=\theta+a,\qquad
 0\le b_q(a)<q^2-1,\qquad a\in\F_q.
\]
The degree of $\theta$ is two, so every parameter is permitted.
Different parameters give different representatives.
If $b_q(a_1)+b_q(a_2)=b_q(a'_1)+b_q(a'_2)$ as integers,
exponentiation gives equal quadratic products at $\theta$.
Their monic leading coefficients cancel, leaving a polynomial
of degree at most one that vanishes at $\theta$.
It is zero, so its two multisets of roots agree. This proves
the Sidon property for the set of representatives.

We next remove all differences at most $2n$ with at most $2n$
deletions. In a Sidon set, a positive difference has at most one
ordered representation: if $u-v=u'-v'>0$, then
$u+v'=u'+v$. Equality of the two unordered pairs implies either
$u=u'$ and $v=v'$, or $u=v$ and $u'=v'$. The second possibility
contradicts positivity. For each integer difference in
$\{1,\ldots,2n\}$ that occurs, mark the larger endpoint of its
unique pair, and delete all marked elements. At most $2n$ elements
are deleted; every pair with a forbidden difference loses an endpoint.
The remaining subset is a Sidon set with the asserted size and gaps.
\end{proof}

The next lemma verifies the recursion in Step~5 for all three
block choices, using only their range and separation properties.

\begin{lemma}\label{lem:iteration}
Suppose the blocks in \eqref{eq:recursion} satisfy Steps~2--3,
with fixed $d$ and $q_k>h\max S_{k-1}$.
Then $U$ is an infinite $B_h$ set containing $A_0$, and
\begin{equation}\label{eq:stage-counts}
 U(n_k)=|S_k|=|S_{k-1}|+|E_{q_k}|,\qquad
 U(M_k-1)=|S_{k-1}|.
\end{equation}
Moreover, $q_k\to\infty$.
\end{lemma}
\begin{proof}
We first verify each finite extension and the infinite union.
At stage $k$, apply Lemma~\ref{lem:extension} with
$n=\max S_{k-1}$, $D=dq_k^h$ and the chosen block $E_{q_k}$.
The integer inequality $q_k>hn$ gives
$M_k\ge(h-1)D+hn+1$.
Induction proves that every $S_k$ is a finite $B_h$ set and that
the new block lies above all old elements.
Each block is nonempty, so $|S_k|\to\infty$ and
$q_k>h\max S_{k-1}\to\infty$.
Any equality of $h$-term sums in $U$ involves finitely many
elements, all contained in one $S_k$.
Its $B_h$ property proves the same property for $U$.

We next check that later blocks do not affect the recorded counts.
All elements of $S_k$ lie below $n_k$, and
$\max S_k\ge M_k$ since the new block is nonempty.
Consequently
\[
 M_{k+1}>q_{k+1}>h\max S_k\ge hM_k>n_k,
 \qquad hM_k-n_k=h(h-2)dq_k^h+(h-1)q_k>0.
\]
Thus $U\cap[1,n_k]=S_k$.
Before $M_k$ only $S_{k-1}$ contributes; the new and old sets
are disjoint. These observations prove \eqref{eq:stage-counts}.
\end{proof}

\section{Density and the main theorem}\label{sec:density}

The common endpoint calculation and selection argument below prove
Theorem~\ref{thm:main}. The final proposition then bounds the count
inside each moment-curve block to obtain its exact upper limit $1/h$.

\begin{lemma}\label{lem:source-density}
Under the hypotheses of Lemma~\ref{lem:iteration}, if
$|E_{q_k}|/q_k\to1$, then the three limits in
\eqref{eq:source-scales} hold.
In particular, $U$ has normalized upper limit at least $(hd)^{-1/h}$
and normalized lower limit zero.
\end{lemma}
\begin{proof}
We first bound the old count relative to the new block size.
By \eqref{eq:parameter-order}, applied with $S=S_{k-1}$,
$n=\max S_{k-1}$ and $q=q_k>hn$,
\begin{equation}\label{eq:old-negligible}
 \frac{|S_{k-1}|}{q_k}
 \le\frac{(h!hn)^{1/h}}{q_k}
 <(h!)^{1/h}q_k^{1/h-1}\longrightarrow0.
\end{equation}
Here $q_k\to\infty$ by Lemma~\ref{lem:iteration}.

We determine the limits by dividing the counts in
\eqref{eq:stage-counts} by the corresponding scales.
The recursion gives
\[
 \frac{n_k^{1/h}}{q_k}\longrightarrow(hd)^{1/h},
 \qquad
 \frac{(M_k-1)^{1/h}}{q_k}\longrightarrow((h-1)d)^{1/h}.
\]
Together with \eqref{eq:old-negligible} and
$|E_{q_k}|/q_k\to1$, these prove the first two limits
in \eqref{eq:source-scales}.
Taking the ratio of $n_k$ and $M_k-1$ in
\eqref{eq:recursion} gives the third.
Nonnegativity then makes the normalized lower limit zero.
\end{proof}

The next lemma selects an exact upper limit while retaining a
prescribed subset with negligible counting function.

\begin{lemma}\label{lem:selection}
Let $U\subseteq\N$ be an infinite $B_h$ set, let $P\subseteq U$
satisfy $P(x)=o(x^{1/h})$, and let
$b=\limsup_{x\to\infty}U(x)/x^{1/h}>0$.
For every $0<c\le b$, there is a $B_h$ set $A$ with
$P\subseteq A\subseteq U$ such that
\begin{equation}\label{eq:selection}
 A(x)=P(x)+\left\lfloor\frac cb\bigl(U(x)-P(x)\bigr)\right\rfloor
     =\frac cbU(x)+o(x^{1/h}).
\end{equation}
If $P$ is finite, the error in the last equality is $O_P(1)$.
The normalized upper limit of $A$ is $c$; if the normalized
lower limit of $U$ is zero, then that of $A$ is also zero.
\end{lemma}
\begin{proof}
We first choose the retained elements outside $P$.
The bound \eqref{eq:counting} makes $b$ finite. The set $U\setminus P$
is infinite, since otherwise $U(x)=P(x)+O(1)=o(x^{1/h})$.
Enumerate its elements in increasing order as $u_1,u_2,\ldots$.
Retain $u_j$ exactly when
$\lfloor cj/b\rfloor>\lfloor c(j-1)/b\rfloor$, and include every
element of $P$. Since $0<c/b\le1$, each increase is either zero
or one. Among the first $l$ elements outside $P$, the number retained
is the telescoping sum $\lfloor cl/b\rfloor$.
Taking $l=U(x)-P(x)$ proves the first equality in \eqref{eq:selection}.

We next determine the counting error and the limiting values.
The first equality gives
$A(x)-(c/b)U(x)=(1-c/b)P(x)+O(1)$.
This is $o(x^{1/h})$, and it is bounded when $P$ is finite.
Division by $x^{1/h}$ now gives the claimed upper and lower limits.
The positive upper limit makes $A$ infinite, and $A\subseteq U$
preserves the $B_h$ property.
\end{proof}

For the improved blocks, let $d_2=1$ and $d_h=h$ for $h\ge3$,
so that $c_h=(hd_h)^{-1/h}$.

\begin{proof}[Proof of Theorem~\ref{thm:main}]
We choose the blocks in \eqref{eq:recursion} with $d=d_h$.
Let $n_0=\max S_0$ and, at stage $k$, choose a prime $q_k>hn_{k-1}$.
For $h\ge3$ use Lemma~\ref{lem:logarithmic}.
For $h=2$ use Lemma~\ref{lem:sidon} with old bound $n_{k-1}$,
and additionally require $q_k>k n_{k-1}$.
The first blocks have $q_k$ elements; the second have between
$q_k-2n_{k-1}$ and $q_k$ elements.
Thus $|E_{q_k}|/q_k\to1$ in both cases.
Their shorter-sum gaps exceed $h n_{k-1}$, which is at least
$h\max S_{k-1}$, so Lemma~\ref{lem:iteration} applies.

We obtain the exact upper limit by selection.
Lemma~\ref{lem:source-density} gives
$b=\limsup U(x)/x^{1/h}\ge c_h$ and lower limit zero.
The endpoint calculation alone does not exclude larger values
inside a block. Apply Lemma~\ref{lem:selection} with $P=A_0$
and $0<c\le c_h\le b$ to obtain the required exact upper limit
and retain the initial set.
\end{proof}

For the moment-curve blocks, the upper limit is exact before any
selection. We recall Proposition~3.2 of our first version, which
establishes this for every admissible choice of primes.

\begin{proposition}\label{prop:moment-density}
With the blocks in \eqref{eq:moment-block} and
$d=h^{h-1}$, the recursion \eqref{eq:recursion} gives an
infinite $B_h$ set $U\supseteq A_0$ satisfying
\[
 \limsup_{x\to\infty}\frac{U(x)}{x^{1/h}}=\frac1h,\qquad
 \liminf_{x\to\infty}\frac{U(x)}{x^{1/h}}=0.
\]
Any primes $q_k>h\max S_{k-1}$ may be used.
\end{proposition}
\begin{proof}
We first determine the endpoint density and the lower limit.
Lemma~\ref{lem:moment} gives $|E_{q_k}|=q_k$.
Lemmas~\ref{lem:iteration} and~\ref{lem:source-density},
with $d=h^{h-1}$, therefore give upper limit at least
$(h^h)^{-1/h}=1/h$ and lower limit zero.

We next bound the count at every point in a new block.
Fix $k$, abbreviate $q=q_k$, and consider $x=M_k+f_q(a)$
with $0\le a<q$. Lemma~\ref{lem:iteration} ensures that no later block
contributes to $U(x)$, and $f_q$ is strictly increasing by
Lemma~\ref{lem:moment}. Hence
$U(x)=|S_{k-1}|+a+1$.
The lower bound in \eqref{eq:moment-range} gives, for $t=a/q$,
\[
 x\ge h^{h-1}q^h(h-1+t),\qquad
 \frac{U(x)}{x^{1/h}}
 \le\frac{t+(|S_{k-1}|+1)/q}
          {h^{(h-1)/h}(h-1+t)^{1/h}}.
\]
For $0\le t\le1$, we have $h-1+t\ge ht$.
Thus the term involving $t$ is at most
$t^{1-1/h}/h\le1/h$, with value zero at $t=0$.
The remaining term is at most
\[
 \frac{(|S_{k-1}|+1)/q}{(h^{h-1}(h-1))^{1/h}}
 \longrightarrow0
\]
by \eqref{eq:old-negligible}, independently of $a$.
Hence the block estimate is uniform.

We finally extend the estimate to all sufficiently large real $x$.
If $u$ is the last element of $U$ at most $x$, then
$U(x)=U(u)$ and $U(x)/x^{1/h}\le U(u)/u^{1/h}$.
As $x\to\infty$, the element $u$ tends to infinity and its block
index does also, because each stage is finite.
The uniform error consequently tends to zero.
This proves the upper bound $1/h$ and hence equality.
\end{proof}

For $h=3$, the moment-curve block and translation are explicitly
\[
 E_q=\{9q^2a+3q[a^2]_q+[a^3]_q:0\le a<q\},
 \qquad M=18q^3+q,\qquad q>3\max S.
\]
The proposition gives upper limit $1/3$ for this construction;
Theorem~\ref{thm:main} permits every value up to $3^{-2/3}$.
Both conclusions concern the full-limit assertion discussed in
the introduction.

\section{Growth between dense scales}\label{sec:growth}

The size of $q_k$ in \eqref{eq:recursion} determines how far the next
block is from the preceding one. For a finite prescribed set,
\eqref{eq:selection} changes the counting function by a fixed
positive factor and a bounded error. Consequently the growth
estimates below pass from the source set to the selected set.
Unless stated otherwise, we use the blocks from the proof of
Theorem~\ref{thm:main}, with $d=d_h$.
The notation $\asymp$ means upper and lower bounds by positive
constant multiples, with fixed parameters permitted in the constants.

\begin{proposition}\label{prop:growth}
Let $h\ge3$, let $A_0\subseteq\N$ be a finite $B_h$ set, and let
$0<c\le c_h$ and $t\ge1$. The set in Theorem~\ref{thm:main} can
be chosen so that
\[
 A(x)\gg x^{1/(h^2t)},\qquad
 \liminf_{x\to\infty}\frac{\log A(x)}{\log x}=\frac1{h^2t},
 \qquad
 \limsup_{x\to\infty}\frac{\log A(x)}{\log x}=\frac1h.
\]
The implied constant may depend on $h,t,c,A_0$.
\end{proposition}
\begin{proof}
We first choose the prime sizes to control the successive endpoints.
Bertrand's theorem permits a prime between
$\lceil hn_{k-1}^t\rceil$ and $2\lceil hn_{k-1}^t\rceil$.
Thus in \eqref{eq:recursion} we can require
$hn_{k-1}^t<q_k<2hn_{k-1}^t+2$.
Since $t\ge1$, the original condition $q_k>hn_{k-1}$ holds.
The endpoint formula and \eqref{eq:source-scales} give
\begin{equation}\label{eq:growth-scales}
 n_k\asymp q_k^h\asymp n_{k-1}^{ht},
 \qquad |S_k|\sim q_k\asymp n_k^{1/h}.
\end{equation}
The constants in these comparisons are independent of $k$.

We next prove the lower bound at every sufficiently large scale.
For $n_{k-1}\le x<n_k$, monotonicity of the count and
\eqref{eq:growth-scales} yield
\[
 U(x)\ge |S_{k-1}|
 \gg n_{k-1}^{1/h}
 \gg n_k^{1/(h^2t)}
 \ge x^{1/(h^2t)}.
\]
Let $b=\limsup U(x)/x^{1/h}$ and apply
Lemma~\ref{lem:selection} with $P=A_0$.
Then $A(x)=(c/b)U(x)+O_{A_0}(1)$, where $c/b>0$ is fixed.
Since $U(x)\to\infty$, the bounded error can be absorbed into a
smaller positive constant in the preceding lower bound.

We identify a sequence on which the lower logarithmic exponent is
attained. Immediately before the $k$th block,
$U(M_k-1)=|S_{k-1}|$. Equations
\eqref{eq:recursion} and \eqref{eq:growth-scales} give
\[
 \log |S_{k-1}|=\frac1h\log n_{k-1}+O(1),\qquad
 \log(M_k-1)=ht\log n_{k-1}+O(1).
\]
Their quotient tends to $1/(h^2t)$.
The bounded-error selection changes the numerator by $O(1)$,
because $|S_{k-1}|\to\infty$.
Together with the all-scale lower bound, this proves the stated
lower logarithmic exponent.

We finally determine the upper logarithmic exponent.
The counting bound \eqref{eq:counting} gives
$\log A(x)/\log x\le1/h+o(1)$.
Since the normalized upper limit of $A$ is $c>0$, along a sequence
$A(x)/x^{1/h}\to c$; on that sequence the logarithmic quotient
tends to $1/h$. This proves the remaining equality.
\end{proof}

For $h=2$, the constants require a separate distinction.
The logarithmic block in Lemma~\ref{lem:logarithmic} has $q-1$
elements in $[0,2q^2)$, so the same proof gives every
$0<c\le1/2$ and exponent $1/(4t)$ for $t\ge1$.
For the full range $0<c\le1/\sqrt2$, use Lemma~\ref{lem:sidon}
and take $t>1$. The choice $q_k\asymp n_{k-1}^t$ then gives
$n_{k-1}/q_k\to0$, as required to make the deletions negligible.
At the endpoint exponent $1/4$, choose instead
$q_k\asymp n_{k-1}\log n_{k-1}$ once $n_{k-1}>1$.
Here
\[
 \log n_k=2\log n_{k-1}
          +2\log\log n_{k-1}+O(1),\qquad
 \log|S_{k-1}|=\tfrac12\log n_{k-1}+O(1).
\]
We obtain an explicit lower bound on the intervals between endpoints.
The same choice of primes gives
\[
 |S_{k-1}|\asymp n_{k-1}^{1/2}
 \asymp\frac{n_k^{1/4}}{\sqrt{\log n_k}}.
\]
Since $x^{1/4}/\sqrt{\log x}$ is increasing for $x>e^2$,
the inequality $U(x)\ge|S_{k-1}|$ for $n_{k-1}\le x<n_k$ gives
$U(x)\gg x^{1/4}/\sqrt{\log x}$.
The fixed positive proportion and bounded error in \eqref{eq:selection}
therefore give $A(x)\gg x^{1/4}/\sqrt{\log x}$.
The values at $M_k-1$ still give lower logarithmic exponent $1/4$.

Primes may also be chosen much farther apart. The next proposition
imposes a slowly growing critical sum and small counts at prescribed
late scales in the same construction. The functions need not be
monotone. For $h=2$, the existence of a Sidon set with
$\sum_{a\in A}a^{-1/2}=\infty$ was proved by
Riblet and Schehr~\cite[Corollary~3.4]{riblet-schehr}. For $h\ge3$,
divergence of $\sum_{a\in A}a^{-1/h}$ already follows from the positive
upper limit in our first version. Here we impose arbitrarily slow
divergence together with a prescribed finite initial set and normalized upper limit.

\begin{proposition}\label{prop:slow-scales}
Let $A_0$ and $c$ satisfy Theorem~\ref{thm:main}. Given positive
functions $g_1,g_2$ tending to infinity and any sequence $y_k\to\infty$,
the set $A$ can satisfy \eqref{eq:limits} and
\[
 \sum_{a\in A}a^{-1/h}=\infty,\qquad
 \sum_{\substack{a\in A\\a\le x}}a^{-1/h}=o(g_1(x)).
\]
It can also have scales $x_k<z_k$, with $x_k\ge y_k$, such that
\[
 \frac{A(x_k)}{g_2(x_k)}\longrightarrow0,\qquad
 \frac{A(z_k)}{z_k^{1/h}}\longrightarrow c,\qquad
 \limsup_{k\to\infty}\frac{z_k}{x_k}\le\frac h{h-1}.
\]
\end{proposition}
\begin{proof}
We first impose the size requirements during the recursion
\eqref{eq:recursion}. Once $S_{k-1}$ is fixed, choose $q_k$ large
enough that
\begin{equation}\label{eq:slow-choices}
 M_k-1\ge\max_{j\le k}y_j,\qquad
 g_2(M_k-1)\ge k|S_{k-1}|,\qquad
 g_1(x)\ge k^2\quad\hbox{for every }x\ge M_k.
\end{equation}
Also keep $q_k>hn_{k-1}$ and, for $h=2$, $q_k>k n_{k-1}$.
These choices are possible because both functions tend to infinity:
each required lower bound holds at every argument beyond a threshold,
and $M_k\to\infty$ as $q_k\to\infty$.
This uses no monotonicity assumption.

We next bound the critical sum of the source set.
Every new block has at most $q_k$ elements, all at least
$M_k\ge q_k^h$. Its contribution to $\sum a^{-1/h}$ is therefore
at most $q_k(q_k^h)^{-1/h}=1$.
For $M_k\le x<M_{k+1}$, only the initial set and the first $k$
blocks contribute, so
\[
 \frac{\sum_{a\in U,\ a\le x}a^{-1/h}}{g_1(x)}
 \le \frac{\sum_{a\in S_0}a^{-1/h}+k}{k^2}\longrightarrow0.
\]
Select $A$ by Lemma~\ref{lem:selection}, with $P=A_0$.
The partial-sum bound passes to $A\subseteq U$, and
\eqref{eq:limits} holds.

We prove divergence from the positive normalized upper limit.
For fixed $y$ and $x>y$, each $a\le x$ satisfies
$x^{-1/h}\le a^{-1/h}$, whence
\[
 \frac{A(x)}{x^{1/h}}
 \le \frac{A(y)}{x^{1/h}}
       +\sum_{\substack{a\in A\\a>y}}a^{-1/h}.
\]
If the series converged, taking the upper limit as $x\to\infty$
and then letting $y\to\infty$ would give $c\le0$.
This contradicts $c>0$ and proves divergence.

We finally locate the requested dense and sparse scales after the
selection. Between consecutive elements of $A$, the count is constant
and its quotient by $x^{1/h}$ decreases. Thus the upper limit $c$
can be approached at elements of $A$.
Choose such elements $z_j$ in blocks with strictly increasing
indices $k_j$. This is possible because each block is finite.
Let $x_j=M_{k_j}-1$. Then \eqref{eq:slow-choices} gives
\[
 x_j\ge y_j,\qquad
 \frac{A(x_j)}{g_2(x_j)}
 \le\frac{|S_{k_j-1}|}{g_2(M_{k_j}-1)}
 \le\frac1{k_j}\longrightarrow0.
\]
Here $k_j\ge j$ explains why the prefix maximum in
\eqref{eq:slow-choices} suffices even if $(y_j)$ is not increasing.
Since $M_{k_j}\le z_j<n_{k_j}$, we have $x_j<z_j$ and
$z_j/x_j\le n_{k_j}/(M_{k_j}-1)$.
Equation~\eqref{eq:source-scales} proves the ratio bound.
\end{proof}

Taking $g_2(x)=\log(x+1)$ gives lower logarithmic exponent zero:
at the scales $x_k$, eventually $A(x_k)\le\log(x_k+1)$,
whereas $\log A(x)/\log x\ge0$ for all sufficiently large $x$.
Together with Proposition~\ref{prop:growth}, this realizes every
lower logarithmic exponent in $[0,1/h^2]$ for $h\ge3$,
while retaining any prescribed $c\in(0,c_h]$.
For comparison, every $B_h$ set satisfies
$\sum_{a\in A}a^{-s}<\infty$ when $s>1/h$: by
\eqref{eq:counting}, the contribution from $[2^j,2^{j+1})$
is $O_h(2^{j(1/h-s)})$, a summable geometric bound.

\begin{corollary}\label{cor:cluster}
For every $B_h$ set satisfying \eqref{eq:limits}, the subsequential
limits of $A(x)/x^{1/h}$ form exactly the interval $[0,c]$.
\end{corollary}
\begin{proof}
We first show that consecutive normalized counts have vanishing
differences. Since $A(n+1)-A(n)\in\{0,1\}$ and
$A(n)=O_h(n^{1/h})$,
\[
 \left|\frac{A(n+1)}{(n+1)^{1/h}}-\frac{A(n)}{n^{1/h}}\right|
 \le (n+1)^{-1/h}
 +O_h(1)\left(1-\left(\frac n{n+1}\right)^{1/h}\right)
 \longrightarrow0.
\]

We next obtain every intermediate limit by a crossing argument.
Fix $0<s<c$. The lower and upper limits in \eqref{eq:limits}
give arbitrarily late indices with values below $s$ and above $s$.
Choose a value below $s$ and a later value above $s$, and take
the first index in between with value at least $s$.
Its distance from $s$ is at most the preceding consecutive
difference, which tends to zero as the starting index tends to
infinity. Thus $s$ is a subsequential limit. The endpoints $0,c$
are subsequential limits by their definitions, and no value outside
$[0,c]$ is permitted by \eqref{eq:limits}.
For real $x$, setting $n=\lfloor x\rfloor$ gives $A(x)=A(n)$
and $n/x\to1$, so the real and integer subsequential limits coincide.
\end{proof}

\section{Ambient restrictions and initial sets}\label{sec:ambient}

Lemma~\ref{lem:extension} remains applicable if a block is replaced
by a nonempty subset or its translation is increased.
The first change preserves the $B_h$ property and all shorter-sum
gaps; the second preserves the required lower bound for the
translation. We use these two facts to place blocks inside a
prescribed set $E\subseteq\N$.
The notation $E_q$ continues to mean a finite block, whereas $E$
denotes the ambient set.

\begin{proof}[Proof of Theorem~\ref{thm:density-one}]
We first find a translate that loses only a small proportion of a
block. Let the current finite set $S$ lie in $E\cap[1,n]$,
fix $0<w<1$, and use the blocks of
Theorem~\ref{thm:main}, with width $d_hq^h$.
Choose the prime $q>hn$ sufficiently large that
\[
 \frac nq\le w,\quad
 q\le(1-w)d_hq^h,\quad
 w d_hq^h\ge1,\quad
 |[1,(h+1)d_hq^h]\setminus E|\le w^2d_hq^h.
\]
The last requirement follows from density one with fixed
$w$. Let $M=(h-1)d_hq^h+q$.
For every integer $t$ between $M$ and
$M+\lfloor w d_hq^h\rfloor$, all of $t+E_q$
lies in $[1,(h+1)d_hq^h]$.

We prove the retained-cardinality estimate by averaging.
For a fixed $e\in E_q$, distinct translations give distinct integers
$t+e$, so at most $w^2d_hq^h$ translations send $e$
outside $E$. Summing this bound over $e\in E_q$, and using at least
$w d_hq^h$ available translations, gives a translation $t$
for which at most $w|E_q|$ elements are lost. Hence
\begin{equation}\label{eq:retained-density-one}
 |(t+E_q)\cap E|\ge(1-w)|E_q|,
 \quad M\le t\le M+w d_hq^h,
 \quad (t+E_q)\cap E\subseteq[1,(h+w)d_hq^h+q).
\end{equation}
This nonempty retained block can be appended to $S$ by
Lemma~\ref{lem:extension}.

We next iterate the retained blocks while controlling their endpoints.
Start from $A_0$, or one point of $E$ if $A_0$ is empty, and choose
$w_k\to0$. At each stage let $n$ be the preceding integer
endpoint, and take the next endpoint to be
$\lceil(h+w_k)d_hq_k^h+q_k\rceil$.
The next translation exceeds this endpoint because its prime exceeds
$hn$. Thus later stages do not change counts at earlier endpoints.
The logarithmic blocks have $q_k$ elements for $h\ge3$;
for $h=2$, the condition $n/q_k\le w_k$ gives
$|E_{q_k}|\ge(1-2w_k)q_k$.
It follows from \eqref{eq:retained-density-one} that each added
block has $(1-o(1))q_k$ elements. Its endpoint divided by $q_k^h$
tends to $hd_h$, so the source upper limit is at least $c_h$.

We finally prove the zero lower limit and select the exact upper limit.
Immediately before each chosen translation $t_k$, only the old set
contributes. By \eqref{eq:counting} and $q_k>hn$ its cardinality is
$O_h(n^{1/h})=o(q_k)$, whereas
$(t_k-1)^{1/h}\gg_h q_k$.
Thus the source has normalized lower limit zero.
Every finite stage is a $B_h$ set by the extension lemma.
Any comparison of two $h$-term sums in the nested union involves
finitely many elements and hence lies in one finite stage, so the
union is a $B_h$ set.
Lemma~\ref{lem:selection}, with $P=A_0$, yields the requested
$c$ while retaining containment in $E$.
\end{proof}

Positive upper density supplies sufficiently many points at a sequence
of endpoints, although it does not make the losses in a short
translation interval negligible. Averaging over a longer interval
gives the following explicit, smaller constant.

\begin{proposition}\label{prop:positive-density}
Let $E\subseteq\N$ have upper density $v>0$, and let
$A_0\subseteq E$ be a finite $B_h$ set. For every
\[
 0<c\le\frac{h c_hv^{1+1/h}}
 {2(h+1)(h+1+v/h)^{1/h}},
\]
there is an infinite $B_h$ set $A$ with $A_0\subseteq A\subseteq E$
satisfying \eqref{eq:limits}.
In particular, $0<c\le c_hv^{1+1/h}/6$ is permitted.
\end{proposition}
\begin{proof}
We first choose a block size relative to a dense endpoint.
Fix the old bound $n$ and
$0<w<hv/(4(h+1))$.
There are arbitrarily large integers $x$ with
$E(x)\ge(v-w)x$.
Bertrand's theorem permits a prime
\begin{equation}\label{eq:positive-prime}
 \left(\frac{v x}{2^h h(h+1)d_h}\right)^{1/h}<q
 \le2\left(\frac{v x}{2^h h(h+1)d_h}\right)^{1/h}.
\end{equation}
Taking such an $x$ large enough ensures
$q>hn$, $n/q\le w$ and $q\le w x$.
The upper bound on $q$ also gives
$hd_hq^h\le v x/(h+1)$.
Use the appropriate block of Theorem~\ref{thm:main} and let
$M=(h-1)d_hq^h+q$.

We next average over all integer translations $M\le t\le x$.
For $e\in E_q\subseteq[0,d_hq^h)$ and
$y\in E\cap[hd_hq^h+q,x]$, the integer $t=y-e$ satisfies
$M\le t\le x$. For fixed $e$, different $y$ give different
translations. Counting these pairs yields
\[
 \begin{aligned}
 \sum_{t=M}^x |(t+E_q)\cap E|
 &\ge |E_q|\bigl(E(x)-hd_hq^h-q\bigr)\\
 &\ge\left(\frac{hv}{h+1}-2w\right)x|E_q|.
 \end{aligned}
\]
There are at most $x$ translations in the sum, so at least one
retains the proportion $hv/(h+1)-2w>0$ of $E_q$.
The retained set lies below $x+d_hq^h$, and its union with the old
set is a $B_h$ set by Lemma~\ref{lem:extension}.

We compute the density supplied by this retained block.
As $w\to0$, the condition $n/q\le w$ makes
$|E_q|/q\to1$, including in the Sidon case.
The lower bound on $q$ in \eqref{eq:positive-prime} and the upper
bound $d_hq^h\le v x/(h(h+1))$ give
\[
 \frac{q}{(x+d_hq^h)^{1/h}}
 \ge
 \frac{c_hv^{1/h}}
 {2(h+1+v/h)^{1/h}}.
\]
Multiplication by the limiting retained proportion
$hv/(h+1)$ gives exactly the constant in the statement.
In particular, the upper limit of the source count is at least
that constant.

We complete the recursion and obtain its lower limit.
Repeat the construction with $w_k\to0$, starting from
$A_0$ or one point of $E$, and use $x+d_hq^h$ as the next old
bound. Later blocks then lie above each recorded endpoint.
Immediately before a new translation $t$, only the old set
contributes, with cardinality $O_h(n^{1/h})=o(q)$.
Since $t\ge(h-1)d_hq^h+q$, the normalized count there tends to
zero. The nested union is a $B_h$ subset of $E$, and
Lemma~\ref{lem:selection} with $P=A_0$ gives every stated $c$.
Finally, $h+1+v/h\le2^h$ for $h\ge2$ and $0<v\le1$,
while $h/(h+1)\ge2/3$.
Therefore the displayed constant is at least
$c_hv^{1+1/h}/6$.
\end{proof}

\begin{corollary}\label{cor:affine-colour}
In a positive arithmetic progression of common difference $m$, every
finite $B_h$ initial set extends to an infinite $B_h$ set satisfying \eqref{eq:limits}
for any $0<c\le c_hm^{-1/h}$.
In every $r$-colouring of $\N$, some colour contains a $B_h$ set
satisfying \eqref{eq:limits} with $c=c_h/r$.
\end{corollary}
\begin{proof}
We first transfer the progression assertion to $\N$.
For the progression $\{a+mn:n\ge0\}$, with $a,m\ge1$,
the map $u\mapsto a+m(u-1)$ is a bijection from $\N$.
An equality of two $h$-term sums in its image cancels the common
translation $h(a-m)$ and the nonzero factor $m$.
Thus this map and its inverse on the progression preserve the
$B_h$ property. Apply Theorem~\ref{thm:main} to the inverse image
of the seed, with upper limit $cm^{1/h}\le c_h$.
The image count is $U((x-a)/m+1)$.
As $((x-a)/m+1)/x\to1/m$, its normalized limits are $c$ and zero.

We next obtain the colouring assertion from a finite partition.
Choose a set $U$ with upper limit $c_h$ and lower limit zero.
Let $U_1,\ldots,U_r$ be its colour parts.
Since $U(x)=\sum_i U_i(x)$,
$c_h\le\sum_i\limsup U_i(x)/x^{1/h}$.
One part therefore has upper limit at least $c_h/r$.
It is infinite and inherits the zero lower limit from $U$.
Lemma~\ref{lem:selection}, with empty prescribed set, selects within
this part a set with upper limit exactly $c_h/r$.
The assertion concerns some colour, not a colour specified in advance.
\end{proof}

An infinite prescribed subset requires room for both the dense blocks
and its later elements. The next sufficient conditions ensure this
room. The assumption that the prescribed union is a $B_h$ set is
necessary for any extension.

\begin{proposition}\label{prop:infinite-seed}
Let $E\subseteq\N$ have density one, let $P=\{b_1<b_2<\cdots\}$
lie in $E$, and let $A_0\subseteq E$ be finite with $A_0\cup P$
a $B_h$ set.
If eventually $b_{k+1}>(4h^2)^h b_k^h$, then for every
$0<c\le h^{-2/h}$ there is an infinite $B_h$ set $A$ with
$A_0\cup P\subseteq A\subseteq E$ satisfying \eqref{eq:limits}.
If instead $b_{k+1}/b_k^h\to\infty$, the full range
$0<c\le c_h$ is permitted.
\end{proposition}
The stronger separation enlarges the range of $c$ only for $h=2$,
from $1/2$ to $1/\sqrt2$; for $h\ge3$ both ranges coincide.
\begin{proof}
We first justify the insertion of an isolated prescribed element.
Let $S$ be a finite nonempty $B_h$ set of positive integers and
let $b>h\max S$. In a comparison of two $h$-term sums from
$S\cup\{b\}$, let $r>s$ be the multiplicities of $b$ on the
two sides. The first sum minus the second is at least
$(r-s)b-(h-s)\max S\ge b-h\max S>0$.
Equal sums must therefore have the same multiplicity of $b$.
Cancelling those copies and padding with the same old element
reduces to the $B_h$ property of $S$. Hence $S\cup\{b\}$
is a $B_h$ set.

We next construct the blocks under the fixed gap bound.
Choose $k_0$ so large that
$b_{k_0}>\max(A_0\cup\{0\})$ and all subsequent size conditions
hold, and start with $A_0\cup\{b_1,\ldots,b_{k_0}\}$.
Suppose the current maximum is $b_k$.
Choose a prime $hb_k<q\le2hb_k$ and use
Lemma~\ref{lem:logarithmic} for every $h\ge2$.
This block has $q+O(1)$ elements in $[0,hq^h)$.
No estimate $b_k/q\to0$ is required here, since no Sidon deletion
is performed.

We place almost all of this block in $E$ with a small additional
translation. Density one and $q\to\infty$ permit
$w_k\to0$ such that
\[
 w_k hq^h\ge1,\qquad
 |[1,(h+1)hq^h]\setminus E|\le w_k^2 hq^h.
\]
For all large $k$, a possible choice is
$w_k=\max\{(hq^h)^{-1/2},
(|[1,(h+1)hq^h]\setminus E|/(hq^h))^{1/2}\}$.
Since $q/(hq^h)\to0$, we also have
$q\le(1-w_k)hq^h$ for all large $k$.
The averaging argument proving \eqref{eq:retained-density-one},
with width $hq^h$, gives a compatible block of $(1-o(1))q$
elements below $(h+w_k)hq^h+q$.
The separate condition $n/q\to0$ from the Sidon case is not used
in this averaging argument.

We check that the next prescribed element can be inserted.
For all sufficiently large $k$, the preceding endpoint is at most
$(h+1)hq^h+q$, and $q\le2hb_k$ gives
\[
 h((h+1)hq^h+q)
 \le 2h^3q^h
 \le (2hq)^h
 \le (4h^2)^h b_k^h
 < b_{k+1}.
\]
The first inequality uses $q\le q^h$ and
$h^2+h+1\le2h^2$.
For the second, $2h^3\le(2h)^h$ is an equality when $h=2$,
and follows from $(2h)^h\ge(2h)^3$ when $h\ge3$.
Thus $b_{k+1}$ exceeds $h$ times the new maximum.
The first paragraph permits its insertion and restores the old
maximum to $b_{k+1}$ for the next stage.

We treat the stronger gap condition using larger intermediate primes.
Choose
\[
 (b_k^h b_{k+1})^{1/(2h)}<q
 \le2(b_k^h b_{k+1})^{1/(2h)}.
\]
Then $q/b_k\to\infty$ and
$q^h/b_{k+1}\le2^h(b_k^h/b_{k+1})^{1/2}\to0$.
For all large $k$ the inequality $q>hb_k$ holds.
Use the blocks of Theorem~\ref{thm:main}, now including the
Sidon deletion when $h=2$.
Its loss is at most $2b_k=o(q)$.
Density-one averaging retains $(1-o(1))q$ elements, with endpoint
$(h+o(1))d_hq^h+q=O_h(q^h)=o(b_{k+1})$.
Increasing $k_0$ if needed ensures that this endpoint is less
than $b_{k+1}/h$ at every stage.
The same singleton argument again inserts $b_{k+1}$.

We determine the normalized limits of these two sources.
Each finite stage is a $B_h$ subset of $E$, so each nested union
is an infinite $B_h$ set containing $A_0\cup P$.
In the first construction the dense endpoints have
$(1-o(1))q$ new elements below $(h^2+o(1))q^h$,
giving upper limit at least $h^{-2/h}$.
In the second the corresponding bound is $(hd_h+o(1))q^h$,
giving upper limit at least $c_h$.
Immediately before each block, the count is that of a finite
$B_h$ set with maximum $b_k$, hence is $O_h(b_k^{1/h})=o(q)$.
The scale is comparable to $q^h$, so both lower limits are zero.

We finally check that selection can retain the entire infinite seed.
Either gap assumption implies $b_{k+1}\ge2b_k^h$ eventually.
For a sufficiently large fixed index $k_0$ with $b_{k_0}>1$,
iteration gives
$\log b_{k_0+j}\ge h^j\log b_{k_0}$.
Thus $b_{k_0+j}\le x$ forces $j=O_h(\log\log(x+3))$,
and $P(x)=O_{h,P}(\log\log(x+3))=o(x^{1/h})$.
Adding the finite set $A_0$ does not change this estimate.
Lemma~\ref{lem:selection}, with prescribed subset $A_0\cup P$,
therefore gives the exact upper limits claimed.
\end{proof}

\begin{proposition}\label{prop:compatible}
Let $E\subseteq\N$ have density one, and prescribe $0<c_i\le c_h$
for each positive integer $i$.
There are pairwise disjoint infinite sets $A_i\subseteq E$, each
satisfying \eqref{eq:limits} with $c=c_i$, whose union is a $B_h$ set.
A finite family of pairwise disjoint finite subsets of $E$ can be
retained in the corresponding sets if its union is a $B_h$ set.
\end{proposition}
\begin{proof}
We first construct the sets together to preserve compatibility.
Begin with the union of the specified seeds, retaining their labels.
If this union is empty, choose one point of $E$ and give it label $1$.
Use the label sequence $1;\ 1,2;\ 1,2,3;\ldots$, in which every
positive integer occurs infinitely often.
At each stage apply the retained-block construction of
Theorem~\ref{thm:density-one} to the whole current finite union,
with $w_k\to0$, and give every element of the new block
that stage's label. All new blocks lie above the previous endpoint,
so the labelled portions $U_i$ are pairwise disjoint.
Every finite union in the recursion is a $B_h$ set; consequently
their infinite union $U$ is also a $B_h$ set.

We next verify the required limits separately for each label.
Immediately before any new block, the total old count is $o(q_k)$,
so $U$ has normalized lower limit zero. Since $U_i\subseteq U$,
every $U_i$ inherits this zero lower limit.
At every endpoint assigned label $i$, its new block alone contributes
$(1-o(1))q_k$ elements below $(hd_h+o(1))q_k^h$.
There are infinitely many such endpoints for each $i$, so
$\limsup U_i(x)/x^{1/h}\ge c_h$ and $U_i$ is infinite.
Apply Lemma~\ref{lem:selection} separately to $U_i$, retaining its
finite seed and selecting upper limit $c_i$.
The resulting sets remain disjoint, and their union is a subset
of $U$, so it retains the $B_h$ property.
\end{proof}

\section{Maximal extensions}\label{sec:maximal}

A $B_h$ subset of $\N$ is maximal if no further positive integer can be
adjoined while preserving the $B_h$ property. For $h=2$, Cilleruelo
and Nathanson~\cite[Theorem~1.3]{cn-perfect} constructed a perfect
difference set $P\subseteq\N$ with $\limsup P(x)/\sqrt{x}\ge1/\sqrt{2}$.
Such a set is maximal Sidon because adjoining $u\notin P$ repeats
the nonzero difference $u-a$ for any $a\in P$. Chen and
Fang~\cite[Theorem~1.1]{chen-fang} proved that for every Sidon set $B$
one can even choose a perfect difference set $P$ with
$|P(x)-B(x/2)|\le\omega(x)$ for any prescribed positive
$\omega(x)\to\infty$.

When only the upper limit is prescribed, the construction in
\cite[Section~3]{cn-perfect} also yields every $0<c\le1/\sqrt2$.
Indeed, put $\alpha=c\sqrt2$. If $l$ is the current maximum, their
prime $p>4l^2$ gives a translated block of $p-O(\sqrt p)$ points in
$p^2+2l+[1,p^2-p]$. Retain $\lfloor\alpha j\rfloor$ of the first
$j$ points for each $j$. At step $k$, adjoin their auxiliary pair
if $k$ is missing from the differences of the actual old set.
The full block and this pair satisfy their Sidon separation bounds,
which are preserved on passing to the retained block. Thus each
positive difference is realized exactly once. A final translation by
one moves their construction, which starts at zero, into $\N$.

Before this final translation, consider $x=p^2+2l+y$ with $0\le y\le p^2$.
The old set contributes at most $l+1=o(p)$ points. The finite Sidon
bound~\cite[Lemma~4]{obryant-I} gives at most $\sqrt y+O(y^{1/4}+1)$
points from the full block up to $x$, hence at most $\alpha\sqrt y+o(p)$
from the retained block, uniformly in $y$. Since
$\sqrt{y/(p^2+y)}\le1/\sqrt2$, the normalized count is at most
$c+o(1)$. At $x=2p^2+2l$ it is $(\alpha p+o(p))/\sqrt{2p^2+2l}
=c+o(1)$. Up to the next block, the auxiliary pair adds at most
two points, so the upper limit is exactly $c$.

The same argument gives the order-two case of Theorem~\ref{thm:maximal}.
Start with the prescribed seed (or one allowed point if it is empty).
Take $\sqrt p>k$ at step $k$ and place every addition beyond the finite
exclusions. Each stage adds $O(1)$ to the critical series. Choose $p$ so
that $g_1(x)\ge k^2$ for $x\ge p^2$ and, at $x_k=p^2+2l$, we have
$x_k\ge y_k$ and the old count is at most $g_2(x_k)/k$. The dense
endpoint is at most $2x_k$. At each stage allow two choices of prime,
with the second block starting beyond all additions of the first choice.
The resulting branches are distinct.

\begin{lemma}\label{lem:exclusion}
Let $S\subseteq\N$ be a finite nonempty $B_h$ set and let
$u\in\N\setminus S$. Above any prescribed bound, one can adjoin
$2h-1$ distinct elements to obtain a $B_h$ set $S'$ such that
$S'\cup\{u\}$ is not a $B_h$ set.
\end{lemma}
\begin{proof}
For a positive integer parameter $t$, define
\[
 a_i=t^{2i},\qquad b_i=t^{2i-1}\quad(1\le i<h),\qquad
 b_h=u+\sum_{i=1}^{h-1}a_i-\sum_{i=1}^{h-1}b_i.
\]
For all sufficiently large $t$, these elements are positive, distinct,
and exceed both the given bound and every element of $S\cup\{u\}$.
Indeed, $b_h$ has leading term $t^{2h-2}$ and differs from
$a_{h-1}$ by a polynomial with leading term $-t^{2h-3}$; all other
distinctions follow from their degrees. The relation
$b_1+\cdots+b_h=u+a_1+\cdots+a_{h-1}$ prevents the insertion of $u$.

Consider two distinct multisets of $h$ labels from $S$ and the new
labels. Their sum difference is a polynomial in $t$. Let $e$ be the
difference of the multiplicities of $b_h$, and let $c_i,d_i$ be those
of $a_i,b_i$ for $i<h$. If the polynomial is identically zero,
comparison of its even and odd powers gives $c_i=-e$ and $d_i=e$.
If $e>0$, one side contains at least $e$ copies of each of the $h$
labels $b_1,\ldots,b_h$. Its length forces $e=1$ and that side to
consist exactly of these labels. The other side then consists of
$a_1,\ldots,a_{h-1}$ and one element $s\in S$; the constant
coefficient gives $u=s$, a contradiction. Exchanging the multisets
excludes $e<0$. If $e=0$, all new multiplicities agree. After their
cancellation, the old multisets agree by the $B_h$ property of $S$,
with identical padding if their common length is less than $h$.
Thus no nontrivial comparison is identically zero. There are finitely
many comparisons and each has finitely many integer roots, so every
sufficiently large integer $t$ outside their union gives the required
$S'$.
\end{proof}

In the following theorem maximality is taken in $\N$.

\begin{theorem}\label{thm:maximal}
Let $A_0\subseteq\N$ be a finite $B_h$ set, and let
$F\subseteq\N\setminus A_0$ be finite. Given $0<c\le c_h$,
positive functions $g_1,g_2$ tending to infinity, and a sequence
$y_k\to\infty$, there are $2^{\aleph_0}$ infinite maximal $B_h$ sets
$A\subseteq\N$ with $A_0\subseteq A$, $A\cap F=\varnothing$,
and \eqref{eq:limits}. Each of these sets also satisfies
\[
 \sum_{a\in A}a^{-1/h}=\infty,\qquad
 \sum_{\substack{a\in A\\a\le x}}a^{-1/h}=o(g_1(x)),
\]
and has integer scales $x_k<z_k$, with $x_k\ge y_k$, such that
\[
 \frac{A(x_k)}{g_2(x_k)}\longrightarrow0,\qquad
 \frac{A(z_k)}{z_k^{1/h}}\longrightarrow c,\qquad
 \limsup_{k\to\infty}\frac{z_k}{x_k}\le\frac h{h-1}.
\]
The functions $g_1,g_2$ need not be monotone.
\end{theorem}
\begin{proof}
Start with $A_0$; if it is empty, start with one integer exceeding
$\max(F\cup\{0\})$. For each $u\in F$, apply
Lemma~\ref{lem:exclusion}, placing all its new elements above $F$
and the current set. We obtain a nonempty finite $B_h$ set $S_0$
containing $A_0$, disjoint from $F$, and preventing the insertion
of every element of $F$.
All elements of $S_0$ are declared forced and retained throughout.
If a $B_h$ superset of $S_0$ contained some $u\in F$, it would
contain $S_0\cup\{u\}$, contradicting the relation just imposed.
Thus every later $B_h$ extension of $S_0$ avoids $F$.

Construct a binary tree of finite sets rooted at $S_0$.
At level $k\ge1$, choose an integer threshold $t_k$ exceeding
$\max S_0$ and, when $k\ge2$, the full candidate endpoint
at every node of level $k-1$. Since each level has finitely
many nodes and $g_1(x)\to\infty$, we may also require
\begin{equation}\label{eq:forced-threshold}
	t_k>\bigl(k(|S_0|+4hk+1)\bigr)^h,
	\qquad g_1(x)\ge k^2\quad\text{for every }x\ge t_k.
\end{equation}
At the first level it also exceeds $\max S_0$.
For the current set $S$, if $k\notin S$ and $S\cup\{k\}$ is a
$B_h$ set, use Lemma~\ref{lem:exclusion} to prevent the insertion
of $k$, with all new terms above $t_k$ and the current set.
If $k$ is already present or already impossible to add, leave $S$
unchanged at this step.

Choose $u>h\max S$ with $u>t_k$. In the first continuation adjoin
$u$; the singleton argument in Proposition~\ref{prop:infinite-seed}
shows that this preserves the $B_h$ property. In the second use
Lemma~\ref{lem:exclusion} to prevent its insertion, placing all
new elements above both $u$ and the current set. Call the elements
added in these two steps forced elements. Each continuation receives
fewer than $4h$ of them at level $k$, all above $t_k$.

For each continuation, let $S$ be its current set and put $n=\max S$.
Choose a prime $q_k>hn$ and the block
from Theorem~\ref{thm:main}, with $d=d_h$. Define $M_k,n_k$
as in \eqref{eq:recursion}, and require
\[
 \frac{|S|}{q_k}\le\frac1k,\qquad
 M_k-1\ge y_k,\qquad g_2(M_k-1)\ge k|S|.
\]
For $h=2$ also require $n/q_k\le1/k$.
All conditions hold for sufficiently large primes, since
$g_2(x)\to\infty$. The chosen block has all distinct shorter sums
separated by more than $hn$, while
$M_k=(h-1)d_hq_k^h+q_k\ge(h-1)d_hq_k^h+hn+1$.
Lemma~\ref{lem:extension}, with $D=d_hq_k^h$, therefore applies to
this current set $S$, including its forced elements.
For $h\ge3$ the block has $q_k$
elements; for $h=2$ it has between $q_k-2n$ and $q_k$ elements.
Consequently $|E_{q_k}|/q_k\to1$.
Process the candidate elements $a$ in increasing order.
If $m$ non-forced elements have already been retained along this
branch, retain $a$ precisely when $m+1\le\lfloor ca^{1/h}\rfloor$.
Forced additions do not change $m$.
Retaining a subset preserves the $B_h$ property. All later additions
must exceed $n_k$, the full candidate endpoint, including the part
containing rejected candidates. There are finitely many nodes at
each level, so these choices define the tree successively.

Fix an infinite branch and let $A$ be the union of its finite sets.
The thresholds tend to infinity, so the corresponding primes do also.
Let $W$ consist of $S_0$ and all its forced elements, and let
$B=A\setminus W$ be the retained non-forced elements.
Through level $k$ there are at most $|S_0|+4hk$ forced elements.
At any forced element of level $k$, \eqref{eq:forced-threshold}
therefore gives $W(x)/x^{1/h}<1/k$. Between consecutive forced
elements this quotient decreases. There is a forced addition at
every level, so $W(x)=o(x^{1/h})$.
The counting threshold gives $B(x)\le\lfloor cx^{1/h}\rfloor$,
and hence $\limsup A(x)/x^{1/h}\le c$.

At a rejected candidate $a$, the preceding count $m$ is at most
$\lfloor ca^{1/h}\rfloor$ by the threshold at earlier candidates,
but $m+1$ exceeds it. Thus $B(a)=\lfloor ca^{1/h}\rfloor$.
Later stages cannot change this value. Infinitely many rejections
give upper limit at least $c$. If only finitely many occur, every
sufficiently late candidate block is retained in full. At such
endpoints,
\begin{equation}\label{eq:maximal-full-block}
 \frac{A(n_k)}{n_k^{1/h}}
 =\frac{|S|/q_k+|E_{q_k}|/q_k}
        {(hd_h+q_k^{1-h})^{1/h}}
 \longrightarrow c_h\ge c.
\end{equation}
In both cases the upper limit equals $c$.
Immediately before each candidate block, $A(M_k-1)=|S|$ and
$(M_k-1)^{1/h}\gg_h q_k$. Our choice $|S|/q_k\le1/k$
therefore gives normalized lower limit zero.

The fewer than $4h$ forced elements of level $k$ contribute at most
$4h/[k(|S_0|+4hk+1)]\le1$ to the critical series by
\eqref{eq:forced-threshold}. A complete candidate block has at most
$q_k$ elements, each at least $q_k^h$, so its contribution is also
at most one. For $x$ from the first new element of level $k$ up to
the first new element of level $k+1$, excluding that last element,
the partial sum is at most $\sum_{a\in S_0}a^{-1/h}+2k$, whereas
$g_1(x)\ge k^2$. These intervals cover all sufficiently large $x$.
Their quotient tends to zero. Positive normalized upper limit forces
divergence by the tail estimate in Proposition~\ref{prop:slow-scales}.

Let $x_k=M_k-1$. Then $x_k\ge y_k$ and
$A(x_k)/g_2(x_k)\le1/k$. If the candidate block has a rejected
element, choose one as $z_k$; the equality for $B(z_k)$ and
$W(z_k)=o(z_k^{1/h})$ give $A(z_k)/z_k^{1/h}\to c$ along these
levels. Otherwise choose $z_k=n_k$.
The calculation in \eqref{eq:maximal-full-block} applies to every
fully retained block, regardless of the decisions at other levels.
Along any infinite sequence of these levels the endpoint ratios
therefore tend to $c_h$.
They occur only finitely often if $c<c_h$, by the upper bound already
proved; if $c=c_h$, they give the required limit as well.
Finally $x_k<z_k\le n_k$ and
$n_k/(M_k-1)\to h/(h-1)$ by their definitions, proving the ratio bound.

Every finite stage is a $B_h$ set, so every branch union is a $B_h$
set by the finite-comparison argument. Each integer $k$ is either
present or impossible to adjoin after level $k$. The finite relation
preventing its insertion is retained forever, so the union is maximal
and still avoids $F$.
At the first differing continuation of two branches, one union
contains the chosen $u$, while the other retains a relation preventing
its insertion. Their unions are different. The binary tree has
$2^{\aleph_0}$ infinite branches, and $\N$ has at most that many
subsets, proving the asserted cardinality.
\end{proof}

Maximality also explains why sparsity alone cannot replace the gap
condition in Proposition~\ref{prop:infinite-seed}. An infinite initial
set can have arbitrarily slow counting growth and still admit no proper
$B_h$ extension. For $h=2$, this also follows from Hall's construction
of an infinite perfect difference set extending a finite Sidon set.
Each pair realizing a missing difference can be placed arbitrarily
far out~\cite[Theorem~3.1]{hall-difference}.

\begin{proposition}\label{prop:sparse-maximal}
Let $A_0\subseteq\N$ be a finite $B_h$ set, and let
$g:[1,\infty)\to(0,\infty)$ be nondecreasing with $g(x)\to\infty$.
There is an infinite maximal $B_h$ set $P\supseteq A_0$ such that
$P(x)\le g(x)$ for all sufficiently large $x$.
In particular, an arbitrary infinite $B_h$ initial set need not admit
an extension with positive normalized upper limit.
\end{proposition}
\begin{proof}
We first construct a maximal set while controlling every new addition.
Start with $A_0$, or with $\{1\}$ if $A_0$ is empty, and process
$k=1,2,\ldots$. For the current finite set $S$, take no action if
$k\in S$ or $S\cup\{k\}$ is not a $B_h$ set. Otherwise choose an
integer $y>\max S$ with $g(y)\ge |S|+2h-1$, and apply
Lemma~\ref{lem:exclusion} with $u=k$, placing all new elements above
$y$. The nested union $P$ is a $B_h$ set. Every omitted integer has
a permanent finite obstruction to its insertion, so $P$ is maximal.
It is infinite, since any finite nonempty $B_h$ set can be enlarged
by an integer greater than $h$ times its maximum.

We check the count on the whole tail, not just at the chosen thresholds.
At the smallest element of a new block, $g$ is already at least the
total cardinality after that block. Until the next block starts,
the count is at most this cardinality and $g$ cannot decrease.
This proves $P(x)\le g(x)$ for every sufficiently large $x$.
For $g(x)=\log(x+2)$, it follows that $P(x)/x^{1/h}\to0$.
Every $B_h$ superset of $P$ equals $P$ by maximality, and therefore
none has positive normalized upper limit.
\end{proof}

The order $h$ is fixed throughout these constructions.
No infinite subset of $\N$ can be a $B_h$ set for unbounded orders.
Indeed, any three of its elements $a<b<c$ satisfy
$(c-b)a+(b-a)c=(c-a)b$.
Both sides have $c-a$ terms, and their multisets are distinct.
Adding the same number of copies of $a$ to both sides gives a
nontrivial equality at every larger order.
The preceding results do not determine the optimal density constant
or necessary and sufficient conditions for retaining an infinite
initial set.

\section*{Declaration of generative AI and AI-assisted technologies}
The authors used a large language model during the preparation of this manuscript
for assistance with mathematical exploration, checking arguments,
organization and language editing. The authors are responsible for
the mathematical statements, proofs, references and final text.

\end{document}